\documentclass[oneside,10pt]{article}          % please do not change
\usepackage[b5paper]{geometry}	    % your paper can be easily printed on a4 or letter paper with enlargenment      
\usepackage{amsfonts,amsmath,latexsym,amssymb} % these packages are required
\usepackage{theorem}                % please use one of this two options for theorems

\usepackage{mathptmx}		    % Journal is printed with poscript fonts: 
\usepackage{jca}	            % journal style, comment only if you find some bug.

\theoremstyle{definition}

\begin{document}

\title[Short title]{Asymptotics of prolate spheroidal wave functions}

% Short title is optional, it will appear in running heads.
% It is necessary only if the title is to long to be used in running heads

\author{T. M. Dunster}

\address{T. M. Dunster, Department of Mathematics and Statistics, 5500 Campanile Drive, San Diego State University, San Diego, CA 92182-7720, USA\\
\email{mdunster@mail.sdsu.edu}}

%  \CorrespondingAuthor{The Name of the corresponding author}

%\dedicated{Dedicated to...}                    % Optional

\date{24.02.2017}                               % Please, write the date of submission


        % AMS-2010 subj class. The list can be found on http://www.ams.org/mathscinet/msc/msc2010.html

%\thanks{This research is supported by ...\par This paper is a lecture that was given at...} 
        % Optional. Only one command thanks is allowed, use \par inside text if you need multiple thanks.       

\begin{abstract}
        Uniform asymptotic approximations are obtained for the
prolate spheroidal wave functions, for both the angular functions
$\operatorname{Ps}_{n}^{m}\left(  {x,\gamma^{2}}\right)  $ (${-1<x<1}$) and
radial functions $Ps_{n}^{m}\left(  {x,\gamma^{2}}\right)  $ (${1<x<\infty}$).
Here $\gamma\rightarrow\infty$, and the results are uniformly valid in the
stated intervals, $m$ and $n$ are integers, with $m$ bounded and $n$
satisfying $0\leq m\leq n\leq2\pi^{-1}\gamma\left(  {1-\delta}\right)  $,
where $\delta\in\left(  {0,1}\right)  $ is fixed. The results are obtained by
an application of certain existing asymptotic solutions of differential
equations, and involve elementary, Bessel, and parabolic cylinder functions.
An asymptotic relationship between the prolate spheroidal equation separation
parameter and the other parameters is also obtained, and error bounds are
available for all approximations.
\end{abstract}

\maketitle

%%%%% END OF TITLE PAGE %%%%%%%%%%%%%%%%%%%%%%%%%%%%%%%%%%%%%%%%%%%%%%

%%%%% BODY OF THE PAPER %%%%%%%%%%%%%%%%%%%%%%%%%%%%%%%%%%%%%%%%%%%%%%
% You should eventually delete (after reading) the rest of the text below %%
 
\numberwithin{equation}{section}
\section{Introduction}

Separation of the wave equation in prolate spheroidal coordinates leads to the
prolate spheroidal wave equation (PSWE)
\begin{equation}
\left(  {1-z^{2}}\right)  \frac{d^{2}y}{dz^{2}}-2z\frac{dy}{dz}+\left(
{\lambda-\frac{\mu^{2}}{1-z^{2}}+\gamma^{2}\left(  {1-z^{2}}\right)  }\right)
y=0, \label{eq1}%
\end{equation}
where $\lambda$ and $\mu$ are separation constants, and $\gamma$ is
proportional to the frequency (see [30] and [40]).

Solutions of \eqref{eq1}, the prolate spheroidal wave functions (PSWFs), are
viewed as depending on the parameters $\mu$ and $\gamma$ from the equation, as
well as an implicitly defined parameter $\nu$ (which describes the behavior of
solutions at infinity). This latter parameter is the so-called characteristic
exponent, and for details see [1, \S 8.1.1].

The parameter $\lambda$ is usually regarded as an eigenvalue admitting an
eigensolution that is bounded at both $z=\pm1$, which is equivalent to $\mu=m$
and $\nu=n$ being integers (see [1] and [19]). Most of the literature focuses
on PSWFs with these parameters being integers, since this is the most useful
case in practical applications. We shall assume this as well throughout this paper.

We consider the important case of $\gamma\rightarrow\infty$ (which
corresponds, for example, to high-frequency scattering in acoustics). In this
case it is known [1, p. 186] that $\lambda\rightarrow-\infty$, and we shall
assume this here. With the exception of \S 5, our results will be uniformly
valid for $m$ bounded, $n$ small or large, and specifically
\begin{equation}
0\leq m\leq n\leq2\pi^{-1}\gamma\left(  {1-\delta}\right)  , \label{eq2}%
\end{equation}
where (here and throughout) $\delta\in\left(  {0,1}\right)  $ is arbitrarily chosen.

Although we will consider the case $z$ complex, our primary concern will be
for $z$ real (denoted by $x)$, and in particular the so-called angular
($-1<x<1$) and the radial ($1<x<\infty$) cases.

We shall apply several existing general asymptotic theories, each of which
contains explicit error bounds. The large range of validity \eqref{eq2}, along
with error bounds, signifies a considerable improvement upon existing results.
We will primarily be concerned with $z$ real, but we shall use complex-valued
argument results as needed to obtain our final results.

One of the principal difficulties in the asymptotic and numerical study of
PSWFs is the determination of the eigenvalues, particularly for large value(s)
of the other parameters. We mention that an extensive the theory of PSWFs with
arbitrary complex parameters $\mu$ and $\nu$ was developed in [19] and [20].

PSWFs were first studied by Niven [24] in heat conduction in spheroidal
bodies, and were subsequently investigated by a number of authors (see [3],
[12], [17], [33], [34], [36]). For certain values of the parameters PSWFs are
eigenfunctions of the finite Fourier transform, and hence these functions play
an important role in signal analysis. These band-limited functions are often
encountered in physics, engineering, statistics; see [4], [32], [38], In [22]
the PSWE is also shown to play a fundamental role in Laplace's tidal
equations. PSWFs also have important applications in fluid dynamics [13],
geophysics and theoretical cosmological models [8], atomic and molecular
physics [15], [28], and biophysics [6].

Despite being studied extensively over the decades there are significant gaps
in the literature on the rigorous analysis of their asymptotic behavior, and
their computation is non-trivial, particularly for large values of the
parameters. The literature contains many asymptotic approximations and
expansions (see [2], [7], [14], [21], [23], [31], [35], [40]), but most are
heuristic, all parameters fixed except $\gamma$, and with little or no error
analysis. For computational techniques see [11], [16], [18], [29], [37], [39].

We mention that in [9] rigorous results (with explicit error bounds) for more
than one large parameter were derived for PSWFs, using a theory of a
coalescing turning point and double pole [5], but not for the parameter range
under consideration in this paper. Specifically, in comparison to the current
paper in which $\lambda<0$, in [9] the case $\lambda>0$ was assumed, which
does not have many of the applications described above.

The PSWE (\ref{eq1}) has regular singularities at $z=\pm1$, each with
exponents $\pm{\frac{1}{2}}m$. When $\gamma=0$ the PSWE degenerates into the
associated Legendre equation (regular singularities at $z=\pm1$ and $z=\infty
$), which for $-1<x<1$ has solutions the Ferrers functions $P_{\nu}^{\mu
}\left(  x\right)  $, and for complex $z$ the associated Legendre functions
$P_{\nu}^{\mu}\left(  z\right)  $.

The significant difference is that if $\gamma\neq0$ the PSWE has an irregular
singularity at infinity. In fact, one can show (from an algebraic form of
Floquet's theorem [1, \S 8.1.1]) that there exists a solution $S_{\nu}%
^{\mu\left(  1\right)  }\left(  {z,\gamma}\right)  $ (in the notation of
[19]), with the property
\begin{equation}
S_{\nu}^{\mu\left(  1\right)  }\left(  {ze^{p\pi i},\gamma}\right)
=e^{p\nu\pi i}S_{\nu}^{\mu\left(  1\right)  }\left(  {z,\gamma}\right)  ,
\label{eq3}%
\end{equation}
for any integer $p$. The LHS of (\ref{eq3}) denotes the branch of the function
after completing $p$ negative half-circuits about $z=\infty$ (equivalently,
$p$ positive half-circuits about $z=\pm1$).

This solution can be expressed as an infinite series involving Bessel
functions of the first kind. Specifically, for integral $\mu$ and $\nu$ we
have from [1, \S 8.3]
\begin{equation}%
\begin{array}
[c]{l}%
S_{n}^{m\left(  1\right)  }\left(  {z,\gamma}\right)  =\left(  {\dfrac{\pi
}{2\gamma z}}\right)  ^{1/2}\dfrac{\left(  {z^{2}-1}\right)  ^{-m/2}z^{m}%
}{A_{n}^{m}\left(  {\gamma^{2}}\right)  }\\
\times\sum\limits_{k=-k^{+ }}^{\infty}{a_{n,k}^{m}\left(  {\gamma^{2}}\right)
J_{n+2k+\left(  {1/2}\right)  }\left(  {\gamma z}\right)  },
\end{array}
\label{eq4}%
\end{equation}
where
\begin{equation}
k^{\pm }=\left\lfloor {{\tfrac{1}{2}}\left(  {n\pm m}\right)  }\right\rfloor .
\label{eq5}%
\end{equation}

We parenthetically note that related solutions $S_{n}^{m\left(  j\right)
}\left(  {z,\gamma}\right)  $ ($j=3,4$) are defined below

In (\ref{eq4}) the coefficients are defined by the three-term recurrence
relation
\begin{equation}%
\begin{array}
[c]{l}%
A_{n,k}^{m}\left(  {\gamma^{2}}\right)  a_{n,k-1}^{m}\left(  {\gamma^{2}%
}\right)  +\left\{  {\lambda_{n}^{m}\left(  {\gamma^{2}}\right)  +B_{n,k}%
^{m}\left(  {\gamma^{2}}\right)  }\right\}  a_{n,k}^{m}\left(  {\gamma^{2}%
}\right) \\
+C_{n,k}^{m}a_{n,k+1}^{m}\left(  {\gamma^{2}}\right)  =0,
\end{array}
\label{eq6}%
\end{equation}
where
\begin{equation}
A_{n,k}^{m}\left(  {\gamma^{2}}\right)  =\frac{\left(  {n-m+2k-1}\right)
\left(  {n-m+2k}\right)  }{\left(  {2n+4k-3}\right)  \left(  {2n+4k-1}\right)
}\gamma^{2}, \label{eq7}%
\end{equation}%
\begin{equation}
B_{n,k}^{m}\left(  {\gamma^{2}}\right)  =\frac{2\left[  {\left(
{n+2k}\right)  \left(  {n+2k+1}\right)  +m^{2}-1}\right]  }{\left(
{2n+4k-1}\right)  \left(  {2n+4k+3}\right)  }\gamma^{2}-\left(  {n+2k}\right)
\left(  {n+2k+1}\right)  , \label{eq8}%
\end{equation}
and
\begin{equation}
C_{n,k}^{m}\left(  {\gamma^{2}}\right)  =\frac{\left(  {n+m+2k+1}\right)
\left(  {n+m+2k+2}\right)  }{\left(  {2n+4k+3}\right)  \left(  {2n+4k+5}%
\right)  }\gamma^{2}. \label{eq9}%
\end{equation}
The normalising constant $A_{n}^{m}\left(  {\gamma^{2}}\right)  $ is defined
by
\begin{equation}
A_{n}^{m}\left(  {\gamma^{2}}\right)  =\sum\limits_{k=-k^{+}}^{\infty
}{\left(  {-1}\right)  ^{k}a_{n,k}^{m}\left(  {\gamma^{2}}\right)  }.
\label{eq10}%
\end{equation}
We also remark that the coefficients $a_{n,k}^{m}\left(  {\gamma^{2}}\right)
$ vanish for $k\leq-1-k^{+}$, where $k^{+}$ is defined by \eqref{eq5}.

From \eqref{eq4}, \eqref{eq10} and the well-known behavior of the Bessel
function at infinity [25, Chap. 12, \S 1.2] we observe that $S_{\nu}%
^{\mu\left(  1\right)  }\left(  {z,\gamma}\right)  $ has the important
property
\begin{equation}
S_{\nu}^{\mu\left(  1\right)  }\left(  {z,\gamma}\right)  =\frac{\cos\left\{
{\gamma z-{\frac{1}{2}}\pi\left(  {\nu+1}\right)  }\right\}  }{\gamma
z}\left\{  {1+{O}\left(  {\frac{1}{z}}\right)  }\right\}  \quad\left(
{z\rightarrow\infty}\right)  , \label{eq11}%
\end{equation}
for $\left\vert {\arg\left(  z\right)  }\right\vert <\pi$. With our assumption
that $\nu=n=$ integer, all solutions of (\eqref{eq1}) are single-valued in the
$z$-plane having a cut along the interval $\left[  {-1,1}\right]  $.

Other fundamental solutions at infinity are given by $S_{n}^{m\left(
j\right)  }\left(  {z,\gamma}\right)  $ ($j=3,4$). These are defined by
\begin{equation}%
\begin{array}
[c]{l}%
S_{n}^{m\left(  j\right)  }\left(  {z,\gamma}\right)  =\left(  {\frac{\displaystyle\pi
}{\displaystyle 2\gamma z}}\right)  ^{1/2}\frac{\left( \displaystyle {z^{2}-1}\right)  ^{-m/2} \displaystyle z^{m}%
}{ \displaystyle A_{n}^{m}\left(  {\gamma^{2}}\right)  }\\
\times\sum\limits_{k=-k^{+}}^{\infty}{a_{n,k}^{m}\left(  {\gamma^{2}}\right)
H_{n+2k+\left(  {1/2}\right)  }^{\left(  {j-2}\right)  }\left(  {\gamma
z}\right)  },
\end{array}
\label{eq12}%
\end{equation}
where $H_{\nu}^{\left(  {1,2}\right)  }\left(  z\right)  $ are the Hankel
functions of the first and second kinds, respectively.

The solutions $S_{n}^{m\left(  j\right)  }\left(  {z,\gamma}\right)  $ have
the fundamental properties
\begin{equation}
S_{n}^{m\left(  3\right)  }\left(  {z,\gamma}\right)  =i^{-n-1}\frac
{e^{i\gamma z}}{\gamma z}\left\{  {1+{O}\left(  {\frac{1}{z}}\right)
}\right\}  \quad\left(  {z\rightarrow\infty}\right)  , \label{eq13}%
\end{equation}
and
\begin{equation}
S_{n}^{m\left(  4\right)  }\left(  {z,\gamma}\right)  =i^{n+1}\frac
{e^{-i\gamma z}}{\gamma z}\left\{  {1+{O}\left(  {\frac{1}{z}}\right)
}\right\}  \quad\left(  {z\rightarrow\infty}\right)  , \label{eq14}%
\end{equation}
for $\left\vert {\arg\left(  z\right)  }\right\vert <\pi$. In particular
$S_{n}^{m\left(  3\right)  }\left(  {z,\gamma}\right)  $ is the unique
solution that is recessive in the upper half plane, and $S_{n}^{m\left(
4\right)  }\left(  {z,\gamma}\right)  $ is the unique solution that is
recessive in the lower half plane.

An important connection formula, which comes directly from the corresponding
one relating the $J$ Bessel function to Hankel functions, is given by
\begin{equation}
S_{n}^{m\left(  1\right)  }\left(  {z,\gamma}\right)  ={\tfrac{1}{2}}\left\{
{S_{n}^{m\left(  3\right)  }\left(  {z,\gamma}\right)  +S_{n}^{m\left(
4\right)  }\left(  {z,\gamma}\right)  }\right\}  . \label{eq15}%
\end{equation}
For $n\geq m\geq0$ and $-1<x<1$ there is a solution $\operatorname{Ps}_{n}%
^{m}\left(  {x,\gamma^{2}}\right)  $ defined in terms of Ferrers functions by
(see [1, \S 8.2])
\begin{equation}
\operatorname{Ps}_{n}^{m}\left(  {x,\gamma^{2}}\right)  =\sum\limits_{k=-k^{-}%
}^{\infty}{\left(  {-1}\right)  ^{k}a_{n,k}^{m}\left(  {\gamma^{2}}\right)
}\operatorname{P}{_{n+2k}^{m}\left(  x\right)  }. \label{eq16}%
\end{equation}
This is the unique solution having the property of being recessive at $x=1$,
and in particular has the property
\begin{equation}
\operatorname{Ps}_{n}^{m}\left(  {x,\gamma^{2}}\right)  =K_{n}^{m}\left(
{\gamma^{2}}\right)  \left(  {1-x}\right)  ^{m/2}\left\{  {1+{O}%
\left(  {1-x}\right)  }\right\}  \quad\left(  {x\rightarrow1^{-}}\right)  ,
\label{eq17}%
\end{equation}
where $K_{n}^{m}\left(  {\gamma^{2}}\right)  $ is a constant given by
\begin{equation}
K_{n}^{m}\left(  {\gamma^{2}}\right)  =\frac{\left(  {-1}\right)  ^{m}%
}{2^{m/2}m!}\sum\limits_{k=-k^{-}}^{\infty}{\left(  {-1}\right)  ^{k}%
\frac{\left(  {n+2k+m}\right)  !}{\left(  {n+2k-m}\right)  !}a_{n,k}%
^{m}\left(  {\gamma^{2}}\right)  }. \label{eq18}%
\end{equation}
From \eqref{eq16} and [27, eq. 14.7.17 ] we remark that $\operatorname{Ps}%
_{n}^{m}\left(  {x,\gamma^{2}}\right)  $ also has the fundamental property of
also being bounded at $x=-1$; this is consequence of the characteristic
exponent $\nu=n$ being an integer. In this case the function satisfies the
normalisation condition
\begin{equation}
\int_{-1}^{1}{\left\{  {\operatorname{Ps}_{n}^{m}\left(  {x,\gamma^{2}%
}\right)  }\right\}  ^{2}dx}=\frac{2\left(  {n+m}\right)  !}{\left(
{2n+1}\right)  \left(  {n-m}\right)  !}. \label{eq19}%
\end{equation}
As mentioned above, the separation constant $\lambda=\lambda_{n}^{m}\left(
{\gamma^{2}}\right)  $ is regarded as an eigenvalue for the case $m$ and $n$
integers, and hence admits the eigensolution ${\operatorname{Ps}_{n}%
^{m}\left(  {x,\gamma^{2}}\right)  }$ that is bounded at $x=\pm1$.

For our purposes the PSWE (\ref{eq1}) therefore takes the form
\begin{equation}
\left(  {1-z^{2}}\right)  \frac{d^{2}y}{dz^{2}}-2z\frac{dy}{dz}+\left\{
{\lambda_{n}^{m}\left(  {\gamma^{2}}\right)  -\frac{m^{2}}{1-z^{2}}+\gamma
^{2}\left(  {1-z^{2}}\right)  }\right\}  y=0. \label{eq20}%
\end{equation}
For $x$ real and lying in $\left(  {1,\infty}\right)  $ (the radial case) we
have the following solution of (\ref{eq20})
\begin{equation}
Ps_{n}^{m}\left(  {x,\gamma^{2}}\right)  =\sum\limits_{k=-k^{-}}^{\infty
}{\left(  {-1}\right)  ^{k}a_{n,k}^{m}\left(  {\gamma^{2}}\right)
P_{n+2k}^{m}\left(  x\right)  }, \label{eq21}%
\end{equation}
which also defines $Ps_{n}^{m}\left(  {z,\gamma^{2}}\right)  $ for complex
$z$; in this case $Ps_{n}^{m}\left(  {z,\gamma^{2}}\right)  $ is entire if $m$
is even, and if $m$ odd $\left(  {1-z^{2}}\right)  ^{1/2}Ps_{n}^{m}\left(
{z,\gamma^{2}}\right)  $ is entire.

Now since $m$ and $n$ are integers, it is a straightforward to show from
(\ref{eq21}) that
\begin{equation}
Ps_{n}^{m}\left(  {ze^{\pi i},\gamma^{2}}\right)  =\left(  {-1}\right)
^{n}Ps_{n}^{m}\left(  {z,\gamma^{2}}\right)  , \label{eq22}%
\end{equation}
which is (unique) property of the Floquet solution $S_{n}^{m\left(  1\right)
}\left(  {z,\gamma}\right)  $. Hence
\begin{equation}
S_{n}^{m\left(  1\right)  }\left(  {z,\gamma}\right)  =\left(  {-1}\right)
^{n}\left(  {n-m}\right)  !V_{n}^{m}\left(  \gamma\right)  Ps_{n}^{m}\left(
{z,\gamma^{2}}\right)  , \label{eq23}%
\end{equation}
for some constant $V_{n}^{m}\left(  \gamma\right)  $, and hence from the known
behavior of $S_{n}^{m\left(  1\right)  }\left(  {z,\gamma}\right)  $ at
infinity
\begin{equation}
Ps_{n}^{m}\left(  {z,\gamma^{2}}\right)  =V_{n}^{m}\left(  \gamma\right)
\frac{\sin\left\{  {\gamma z-{\frac{1}{2}}\pi n}\right\}  }{\gamma z}\left\{
{1+{O}\left(  {\frac{1}{z}}\right)  }\right\}  \quad\left(
{z\rightarrow\infty}\right)  . \label{eq24}%
\end{equation}
An explicit expression, in terms of $a_{n,k}^{m}\left(  {\gamma^{2}}\right)
$, for the constant $V_{n}^{m}\left(  \gamma\right)  $ can be obtained from
\eqref{eq11}, \eqref{eq21}, \eqref{eq23}, and letting $z\rightarrow\infty$.

From [1, p. 171] we also note the recessive behavior
\begin{equation}
Ps_{n}^{m}\left(  {z,\gamma^{2}}\right)  =K_{n}^{m}\left(  {\gamma^{2}%
}\right)  \left(  {z-1}\right)  ^{m/2}\left\{  {1+{O}\left(
{z-1}\right)  }\right\}  \quad\left(  {z\rightarrow1}\right)  , \label{eq25}%
\end{equation}
where $K_{n}^{m}\left(  {\gamma^{2}}\right)  $ is given by (\ref{eq18}).

The plan of the paper is as follows. In \S 2 we obtain Liouville-Green
approximations for $S_{n}^{m\left(  j\right)  }\left(  {z,\gamma}\right)  $
($j=3,4)$ where $z$ is complex, and use these to obtain an asymptotic
approximation for the radial PSWF $Ps_{n}^{m}\left(  {x,\gamma^{2}}\right)  $
which is uniformly valid in the interval $1+\delta\leq x<\infty$. In \S 3 the
approximation for $Ps_{n}^{m}\left(  {x,\gamma^{2}}\right)  $ is extended to
$1<x<\infty$ by applying the theory of differential equations having a simple
pole, which involves the Bessel function of the first kind. Also in this
section an asymptotic relationship involving $\lambda_{n}^{m}\left(
{\gamma^{2}}\right)  _{\,}$and the parameters $m$ and $n$ is obtained, by
matching the Liouville-Green and Bessel function approximations at infinity.

In \S 4 the angular PSWF $\operatorname{Ps}_{n}^{m}\left(  {x,\gamma^{2}%
}\right)  $ is approximated, with the intervals $1-\delta_{0}\leq x<1$ and
$0\leq x\leq1-\delta_{0}$ considered separately (for some positive constant
$\delta_{0})$. In the former interval the asymptotic approximation involves
the modified Bessel function of the first kind, and in the latter interval the
parabolic cylinder function is used. In \S 5 the approximation involving the
parabolic cylinder function is simplified under the assumption $n$ being
bounded. Finally, in \S 6 we summarise the main results of the paper.

\section{Liouville-Green asymptotics: the radial case}

Making the transformation $w=\left(  {z^{2}-1}\right)  ^{1/2}y$ in
(\ref{eq20}) we remove the first derivative to obtain
\begin{equation}
\frac{d^{2}w}{dz^{2}}=\left\{  {-\gamma^{2}+\frac{\lambda_{n}^{m}\left(
{\gamma^{2}}\right)  }{z^{2}-1}+\frac{m^{2}-1}{\left(  {z^{2}-1}\right)  ^{2}%
}}\right\}  w. \label{eq26}%
\end{equation}
Now, from [1, p. 186] it is known that for large $\gamma$, with $m$ and $n$
bounded, that
\begin{equation}
\lambda_{n}^{m}\left(  {\gamma^{2}}\right)  =-\gamma^{2}+2\left(
{n-m+{\tfrac{1}{2}}}\right)  \gamma+{O}\left(  1\right)  .
\label{eq27}%
\end{equation}
With this in mind we define a parameter $\sigma$ by
\begin{equation}
\lambda_{n}^{m}\left(  {\gamma^{2}}\right)  =-\gamma^{2}\left(  {1-\sigma^{2}%
}\right)  , \label{eq28}%
\end{equation}
and throughout we shall assume that
\begin{equation}
0\leq\sigma=\sqrt{1+\gamma^{-2}\lambda_{n}^{m}\left(  {\gamma^{2}}\right)
}\leq\sigma_{0}<1, \label{eq29}%
\end{equation}
where $\sigma_{0}$ is an arbitrary positive constant.

Next, from \eqref{eq28} we can express (\ref{eq26}) in the form
\begin{equation}
\label{eq30}\frac{d^{2}w}{dz^{2}}=\left[  {\gamma^{2}f\left(  {\sigma,z}
\right)  +g\left(  z \right)  } \right]  w,
\end{equation}
where
\begin{equation}
\label{eq31}f\left(  {\sigma,z} \right)  =\frac{\sigma^{2}-z^{2}}{z^{2}%
-1},\quad g\left(  z \right)  =\frac{m^{2}-1}{\left(  {z^{2}-1} \right)  ^{2}%
}.
\end{equation}
We observe for large $\gamma$ the differential equation has turning points at
$z=\pm\sigma$, and on account of our assumption (\ref{eq29}) these turning
points lie in the interval $\left(  {-1,1} \right)  $, they may coalesce with
one another at $z=0$, but are bounded away from the poles $z=\pm1$.

We shall construct Liouville-Green approximations for $Ps_{n}^{m} \left(
{z,\gamma^{2}} \right)  $, using the theory of [25, Chap. 10]. To this end, we
introduce a new independent variable
\begin{equation}
\label{eq32}\xi=\int_{1}^{z} {\left\{  {-f\left(  {\sigma,t} \right)  }
\right\}  ^{1/2}dt} =\int_{1}^{z} {\left(  {\frac{t^{2}-\sigma^{2}}{t^{2}-1}}
\right)  ^{1/2}dt} .
\end{equation}
Branch cuts are suitably chosen so that $0\le\xi<\infty$ for $1\le z<\infty$.

The RHS of (\ref{eq32}) can be expressed in terms of the elliptic integral of
the second kind [27, eq. 19.2.5]
\begin{equation}
E\left(  {a;b}\right)  =\int_{0}^{a}{\left(  {\frac{1-b^{2}t^{2}}{1-t^{2}}%
}\right)  ^{1/2}dt}=b\int_{0}^{a}{\left(  {\frac{b^{-2}-t^{2}}{1-t^{2}}%
}\right)  ^{1/2}dt}. \label{eq33}%
\end{equation}
Here $b=\sigma^{-1}>1$, and the branches of the square roots are such that
integrand is positive for $0\leq t<b^{-1}$ and negative for $1<t<\infty$, and
continuous elsewhere in the complex $t$-plane having a cut along the interval
$\left[  {b^{-1},1}\right]  $. We thus have
\begin{equation}
\xi=\sigma E\left(  {z;\sigma^{-1}}\right)  -\sigma E\left(  {1;\sigma^{-1}%
}\right)  . \label{eq34}%
\end{equation}
Then with the new dependent variable $W=\left\{  {-f}\right\}  ^{1/4}w$ we
obtain
\begin{equation}
\frac{d^{2}W}{d\xi^{2}}=\left[  {-\gamma^{2}+\psi\left(  \xi\right)  }\right]
W, \label{eq35}%
\end{equation}
where
\begin{equation}
\psi\left(  \xi\right)  =\frac{m^{2}-1}{\left(  {z^{2}-1}\right)  \left(
{z^{2}-\sigma^{2}}\right)  }+\frac{\left(  {1-\sigma^{2}}\right)  \left(
{6z^{4}-\left(  {3+\sigma^{2}}\right)  z^{2}-2\sigma^{2}}\right)  }{4\left(
{z^{2}-1}\right)  \left(  {z^{2}-\sigma^{2}}\right)  ^{3}}. \label{eq36}%
\end{equation}
We observe that $\psi\left(  \xi\right)  ={O}\left(  {\xi^{-2}%
}\right)  $ as $\xi\rightarrow\infty$, but is unbounded at the singularities
$z=\pm1$, and also at the turning points $z=\pm\sigma$.

From the definition of $\xi$ we find that
\begin{equation}
\xi=z-J\left(  \sigma\right)  +{O}\left(  {z^{-1}}\right)
\quad\left(  {z\rightarrow\infty}\right)  , \label{eq37}%
\end{equation}
where
\begin{equation}
J\left(  \sigma\right)  =1-\int_{1}^{\infty}{\left[  {\left(  {\frac
{t^{2}-\sigma^{2}}{t^{2}-1}}\right)  ^{1/2}-1}\right]  dt}. \label{eq38}%
\end{equation}
Note $J\left(  0\right)  =0$ and $J\left(  1\right)  =1$. Now by Cauchy's
theorem
\begin{equation}
0=\operatorname{Re}\int_{-\infty}^{\infty}{\left[  {\left(  {\frac
{t^{2}-\sigma^{2}}{t^{2}-1}}\right)  ^{1/2}-1}\right]  dt}=2\operatorname{Re}%
\int_{0}^{\infty}{\left[  {\left(  {\frac{t^{2}-\sigma^{2}}{t^{2}-1}}\right)
^{1/2}-1}\right]  dt}. \label{eq39}%
\end{equation}
Hence
\begin{equation}
\int_{1}^{\infty}{\left[  {\left(  {\frac{t^{2}-\sigma^{2}}{t^{2}-1}}\right)
^{1/2}-1}\right]  dt}=-\operatorname{Re}\int_{0}^{1}{\left[  {\left(
{\frac{t^{2}-\sigma^{2}}{t^{2}-1}}\right)  ^{1/2}-1}\right]  dt}, \label{eq40}%
\end{equation}
and consequently from (\ref{eq38})
\begin{equation}
J\left(  \sigma\right)  =1+\operatorname{Re}\int_{0}^{1}{\left[  {\left(
{\frac{t^{2}-\sigma^{2}}{t^{2}-1}}\right)  ^{1/2}-1}\right]  dt}=\int
_{0}^{\sigma}{\left(  {\frac{\sigma^{2}-t^{2}}{1-t^{2}}}\right)  ^{1/2}dt};
\label{eq41}%
\end{equation}
i.e.
\begin{equation}
J\left(  \sigma\right)  =\sigma E\left(  {\sigma;\sigma^{-1}}\right)  ,
\label{eq42}%
\end{equation}
for $\sigma>0$, in which $E$ is the Elliptic integral of the second kind given
by (\ref{eq33}). Thus
\begin{equation}
\xi=z-\sigma E\left(  {\sigma;\sigma^{-1}}\right)  +{O}\left(
{z^{-1}}\right)  \quad\left(  {z\rightarrow\infty}\right)  . \label{eq43}%
\end{equation}
We now apply Theorem 3.1 of [26], with $u$ replaced by $\gamma$, and with
$\xi$ replaced by $i\xi$. Then, by matching solutions that are recessive at
$z=\pm i\infty$, we have from \eqref{eq13}, (\ref{eq14}) and (\ref{eq43})
\begin{equation}%
\begin{array}
[c]{l}%
S_{n}^{m\left(  3\right)  }\left(  {z,\gamma}\right)  =i^{-1-n}\gamma
^{-1}\left[  {\left(  {z^{2}-1}\right)  \left(  {z^{2}-\sigma^{2}}\right)
}\right]  ^{-1/4}e^{i\gamma J\left(  \sigma\right)  }\\
\times\left[  {e^{i\gamma\xi}\sum\limits_{s=0}^{p-1}{\left(  {-i}\right)
^{s}\frac{\displaystyle A_{s}\left(  \xi\right)  }{\displaystyle\gamma^{s}}}+\varepsilon_{p,1}\left(
{\gamma,\xi}\right)  }\right]  ,
\end{array}
\label{eq44}%
\end{equation}
and
\begin{equation}%
\begin{array}
[c]{l}%
S_{n}^{m\left(  4\right)  }\left(  {z,\gamma}\right)  =i^{1+n}\gamma
^{-1}\left[  {\left(  {z^{2}-1}\right)  \left(  {z^{2}-\sigma^{2}}\right)
}\right]  ^{-1/4}e^{-i\gamma J\left(  \sigma\right)  }\\
\times\left[  {e^{-i\gamma\xi}\sum\limits_{s=0}^{p-1}{i^{s}\frac{\displaystyle A_{s}\left(
\xi\right)  }{\displaystyle \gamma^{s}}}+\varepsilon_{p,2}\left(  {\gamma,\xi}\right)
}\right]  .
\end{array}
\label{eq45}%
\end{equation}

The error terms $\varepsilon_{p,j}\left(  {\gamma,\xi}\right)  $ ($j=1,2$) are
bounded by Olver's theorem, and are ${O}\left(  {\gamma^{-p}}\right)
$ in unbounded domains containing the real interval $1+\delta\leq z<\infty$
($\delta>0$). Here the coefficients are defined recursively by $A_{0}\left(
\xi\right)  =1$ and
\begin{equation}
A_{s+1}\left(  \xi\right)  =-{\tfrac{1}{2}A}_{s}^{\prime}\left(  \xi\right)
+{\tfrac{1}{2}}\int{\psi\left(  \xi\right)  A_{s}\left(  \xi\right)  d\xi
}\quad\left(  {s=0,1,2,\cdots}\right)  . \label{eq46}%
\end{equation}
Thus, from (\ref{eq15}), (\ref{eq23}), (\ref{eq42}), (\ref{eq44}) and
(\ref{eq45}), we obtain the desired Liouville-Green expansion for $Ps_{n}%
^{m}\left(  {x,\gamma^{2}}\right)  $. In particular, to leading order, we
have
\begin{equation}
Ps_{n}^{m}\left(  {x,\gamma^{2}}\right)  =\frac{\left(  {-1}\right)  ^{n}%
\sin\left(  {\gamma\xi+\gamma\sigma E\left(  {\sigma;\sigma^{-1}}\right)
-{\frac{1}{2}}n\pi}\right)  +{O}\left(  {\gamma^{-1}}\right)  }%
{\gamma\left(  {n-m}\right)  !V_{n}^{m}\left(  \gamma\right)  \left[  {\left(
{x^{2}-1}\right)  \left(  {x^{2}-\sigma^{2}}\right)  }\right]  ^{1/4}},
\label{eq47}%
\end{equation}
as $\gamma\rightarrow\infty$, uniformly for $1+\delta\leq x<\infty$. In order
for this approximation to be practicable, one requires an asymptotic
approximation for $\lambda_{n}^{m}\left(  {\gamma^{2}}\right)  $ as
$\gamma\rightarrow\infty$, and we shall discuss this in the next section. We
also remark that (\ref{eq47}) breaks down at the simple pole $x=1$, and in the
next section we obtain asymptotic approximations that are valid at this pole.

\section{Bessel function approximations: the radial case}

We now obtain approximations valid at the simple pole of $f\left(  {\sigma
,z}\right)  $ at $z=1$, using the asymptotic theory of [25, Chap. 12]. We
consider $z=x$ real and positive. The appropriate Liouville transformation is
now given by
\begin{equation}
\eta=\xi^{2}=\left[  {\int_{1}^{x}{\left\{  {-f\left(  {\sigma,t}\right)
}\right\}  ^{1/2}dt}}\right]  ^{2}, \label{eq48}%
\end{equation}
along with
\begin{equation}
\hat{{W}}=\left\{  {\frac{\eta\left(  {x^{2}-\sigma^{2}}\right)  }{x^{2}-1}%
}\right\}  ^{1/4}w, \label{eq49}%
\end{equation}
which yields the new equation
\begin{equation}
\frac{d^{2}\hat{{W}}}{d\eta^{2}}=\left[  {-\frac{\gamma^{2}}{4\eta}%
+\frac{m^{2}-1}{4\eta^{2}}+\frac{\hat{{\psi}}\left(  \eta\right)  }{\eta}%
}\right]  \hat{{W}}. \label{eq50}%
\end{equation}
Here
\begin{equation}%
\begin{array}
[c]{l}%
\hat{{\psi}}\left(  \eta\right)  =\dfrac{1-4m^{2}}{16\eta}+\dfrac{m^{2}%
-1}{4\left(  {x^{2}-1}\right)  \left(  {x^{2}-\sigma^{2}}\right)  }\\
+\dfrac{\left(  {1-\sigma^{2}}\right)  \left(  {6x^{4}-\left(  {3+\sigma^{2}%
}\right)  x^{2}-2\sigma^{2}}\right)  }{16\left(  {x^{2}-1}\right)  \left(
{x^{2}-\sigma^{2}}\right)  ^{3}}.
\end{array}
\label{eq51}%
\end{equation}
This has the same main features of (\ref{eq30}), namely a simple pole for the
dominant term (for large $\gamma$) and a double pole in another term. We note
that $x=1$ corresponds to $\eta=0$.

The difference here is that non-dominant term $\hat{{\psi}}\left(
\eta\right)  $ is now analytic at $\eta=0$, i.e. $x=1$. Neglecting $\hat
{{\psi}}\left(  \eta\right)  $ in \eqref{eq50} gives an equation solvable in
terms of Bessel functions. We then find (by matching recessive solutions at
$x=1$) and applying theorem 4.1 of [25, Chap. 12] (with $u$ replaced by
$\gamma$ and $\zeta$ replaced by $\eta$)
\begin{equation}%
\begin{array}
[c]{l}%
Ps_{n}^{m}\left(  {x,\gamma^{2}}\right)  =c_{n}^{m}\left(  \gamma\right)
\left\{  {\dfrac{\eta}{\left(  {x^{2}-1}\right)  \left(  {x^{2}-\sigma^{2}%
}\right)  }}\right\}  ^{1/4}\\
\times\left[  {J_{m}\left(  {\gamma\eta^{1/2}}\right)  +{O}\left(
{\gamma^{-1}}\right)  \operatorname{env}J_{m}\left(  {\gamma\eta^{1/2}%
}\right)  }\right]  ,
\end{array}
\label{eq52}%
\end{equation}
as $\gamma\rightarrow\infty$, uniformly for $1<x<\infty$. Olver's theorem
provides an asymptotic expansion in inverse powers of $\gamma$, but we present
just the leading term here. The so-called envelope env of the $J$ Bessel
function is defined by [27, \S 2.8(iv)].

The constant of proportionality $c_{n}^{m}\left(  \gamma\right)  $ can be
found by comparing both sides of (\ref{eq52}) as $x\rightarrow1$
($\eta\rightarrow0$). Using
\begin{equation}
\eta=2\left(  {1-\sigma^{2}}\right)  \left(  {x-1}\right)  +{O}%
\left\{  {\left(  {x-1}\right)  ^{2}}\right\}  \quad\left(  {x\rightarrow
1}\right)  , \label{eq53}%
\end{equation}
along with (\ref{eq25}), (\ref{eq28}) and the behavior of the $J$ Bessel
function at the origin (e.g. [25, Chap. 12, \S 1]) we find that
\begin{equation}
c_{n}^{m}\left(  \gamma\right)  =\left(  {-\frac{2}{\lambda_{n}^{m}\left(
{\gamma^{2}}\right)  }}\right)  ^{m/2}m!K_{n}^{m}\left(  {\gamma^{2}}\right)
. \label{eq54}%
\end{equation}
An asymptotic approximation for this constant is given by (\ref{eq105}) below.

Next, from the well-known behavior of the $J$ Bessel function at infinity
(e.g. see [25, Chap. 12, \S 1]), we find from \eqref{eq52} that
\begin{equation}%
\begin{array}
[c]{l}%
Ps_{n}^{m}\left(  {x,\gamma^{2}}\right)  \sim\text{constant}\times\left\{
{\left(  {x^{2}-1}\right)  \left(  {x^{2}-\sigma^{2}}\right)  }\right\}
^{-1/4}\\
\times\left\{  {\cos\left(  {\gamma\xi-\frac{1}{2}m\pi-\frac{1}{4}\pi}\right)
+{O}\left(  \xi^{-1}\right)  }\right\}  \quad\left(  {\eta=\xi
^{2}\rightarrow\infty}\right)  .
\end{array}
\label{eq55}%
\end{equation}
However, from (\ref{eq47}) we observe an alternative expression of the
behavior of this function. On comparing both, we deduce that
\begin{equation}
\gamma\sigma E\left(  {\sigma;\sigma^{-1}}\right)  =\left(  {2N+{\tfrac{1}{2}%
}n-{\tfrac{1}{2}}m+{\tfrac{1}{4}}}\right)  \pi+{O}\left(  {\gamma
^{-1}}\right)  , \label{eq56}%
\end{equation}
for some integer $N$, which we show is zero. Now, from (\ref{eq27}) and
(\ref{eq28}) we have for \textit{fixed} $m$ and $n$ that $\sigma
={O}\left(  {\gamma^{-1/2}}\right)  $ as $\gamma\rightarrow\infty$,
and more precisely,
\begin{equation}
\sigma^{2}=2\left(  {n-m+{\tfrac{1}{2}}}\right)  \gamma^{-1}+{O}%
\left(  {\gamma^{-2}}\right)  . \label{eq57}%
\end{equation}
Thus in this case, using (\ref{eq41}) and (\ref{eq42}) in the LHS of
(\ref{eq56}), we have
\begin{equation}
{\tfrac{1}{4}}\pi\gamma\sigma^{2}+{O}\left(  {\gamma^{-1}\sigma^{4}%
}\right)  =\left(  {2N+{\tfrac{1}{2}}n-{\tfrac{1}{2}}m+{\tfrac{1}{4}}}\right)
\pi+{O}\left(  {\gamma^{-1}}\right)  . \label{eq58}%
\end{equation}
Inserting (\ref{eq57}) into (\ref{eq58}) we deduce that $N=0$, at least for
fixed $m$ and $n$; a continuity argument removes this restriction.

It is possible to extend (\ref{eq47}) and (\ref{eq55}) to asymptotic
expansions, and consequently from (\ref{eq41}) and (\ref{eq42}) we arrive at
\begin{equation}
\label{eq59}\gamma\int_{0}^{\sigma}{\left(  {\frac{\sigma^{2}-t^{2}}{1-t^{2}}}
\right)  ^{1/2}dt} \sim\frac{1}{2}\left(  {n-m+\frac{1}{2}} \right)  \pi
+\sum\limits_{s=0}^{\infty}{\frac{\kappa_{s} }{\gamma^{2s+1}}} ,
\end{equation}
for constants $\kappa_{s} $ which can be determined in terms of the
coefficients appearing in (\ref{eq44}) and (\ref{eq45}). From (\ref{eq29}) we
can invert this expansion to provide a means of computing the eigenvalue
$\lambda=\lambda_{n}^{m} \left(  {\gamma^{2}} \right)  $ asymptotically in
terms of $m$ and $n$ as $\gamma\to\infty$.

Now, the elliptic integral on the LHS of (\ref{eq59}) monotonically increases
from 0 to 1 as $\sigma$ increases from 0 to 1. We therefore see that the
condition (\ref{eq29}) (along with our assumption that $m$ is bounded) is
equivalent to (\ref{eq2}).

\section{Bessel and parabolic cylinder function approximations: the angular
case}

Recall $\operatorname{Ps}_{n}^{m}\left(  {x,\gamma^{2}}\right)  $ is the
unique solution with the property of being recessive at $x=\pm1$. It is also
uniquely determined by the property
\begin{equation}
\operatorname{Ps}_{n}^{m}\left(  {-x,\gamma^{2}}\right)  =\left(  {-1}\right)
^{m+n}\operatorname{Ps}_{n}^{m}\left(  {x,\gamma^{2}}\right)  . \label{eq60}%
\end{equation}
Thus, it suffices to approximate $\operatorname{Ps}_{n}^{m}\left(
{x,\gamma^{2}}\right)  $ in the interval $0\leq x<1$. We consider the
subintervals $0\leq x\leq1-\delta_{0}$ and $1-\delta_{0}\leq x<1$ separately,
where $\delta_{0}\in\left(  {0,1-\sigma_{0}}\right)  $ is arbitrary; recall
that $\sigma_{0}$ is defined by (\ref{eq29}). The significance of this choice
is that the turning point $x=\sigma$ is bounded away from the interval
$\left[  {1-\delta_{0},1}\right]  $.

For $1-\delta_{0}\leq x<1$ we apply theorem 3.1 of [25, Chap. 12]. It can then
be shown by utilising the recessive behavior at $x=1$ ($\eta=0$) that
\begin{equation}%
\begin{array}
[c]{l}%
\operatorname{Ps}_{n}^{m}\left(  {x,\gamma^{2}}\right)  =c_{n}^{m}\left(
\gamma\right)  \left\{  {\dfrac{\left\vert \eta\right\vert }{\left(  {1-x^{2}%
}\right)  \left(  {x^{2}-\sigma^{2}}\right)  }}\right\}  ^{1/4}\\
\times I_{m}\left(  {\gamma\left\vert \eta\right\vert ^{1/2}}\right)  \left[
{1+{O}\left(  {\dfrac{1}{\gamma}}\right)  }\right]  ,
\end{array}
\label{eq61}%
\end{equation}
where $I_{m}\left(  x\right)  $ is the modified Bessel function, and
$c_{n}^{m}\left(  \gamma\right)  _{\,}$is given by (\ref{eq54}). Expansions
and error bounds are obtainable from Olver's theorem.

The interval $0\leq x\leq1-\delta_{0}$ is less straightforward. From
(\ref{eq30}) and (\ref{eq31}) we observe that equation has the turning point
$x=\sigma$ in this interval, and this coalesces with the other turning point
$x=-\sigma$ when $\sigma\rightarrow0$. The appropriate asymptotic theory for
this situation is provided by [26], and from eq. (\ref{eq28}) of this
reference the appropriate transformation is given by
\begin{equation}
\frac{d\zeta}{dx}=\left(  {\frac{\sigma^{2}-x^{2}}{\left(  {1-x^{2}}\right)
\left(  {\alpha^{2}-\zeta^{2}}\right)  }}\right)  ^{1/2}. \label{eq62}%
\end{equation}
Upon integration, this yields the implicit relationship
\begin{equation}
\int_{-\alpha}^{\zeta}{\left(  {\alpha^{2}-\tau^{2}}\right)  ^{1/2}d\tau}%
=\int_{-\sigma}^{x}{\left\{  {-f\left(  {\sigma,t}\right)  }\right\}
^{1/2}dt}=\int_{-\sigma}^{x}{\left(  {\frac{\sigma^{2}-t^{2}}{1-t^{2}}%
}\right)  ^{1/2}dt}. \label{eq63}%
\end{equation}
The lower limits are selected to ensure that the turning point $x=-\sigma$ is
mapped to a new turning point at $\zeta=-\alpha$ (see \eqref{eq69} below).
From [26, eq. (\ref{eq30})] we find that $\alpha$ is given by
\begin{equation}
\alpha^{2}=\dfrac{2}{\pi}\int_{-\sigma}^{\sigma}{\left(  {\frac{\sigma
^{2}-t^{2}}{1-t^{2}}}\right)  ^{1/2}dt}=\frac{4}{\pi}J\left(  \sigma\right)  ,
\label{eq64}%
\end{equation}
which ensures that the original turning point $x=\sigma$ is mapped to the
turning point at $\zeta=\alpha$ in the transformed equation.

By symmetry $x=0$ is mapped to $\zeta=0$, and so the lower limits in the
integrals of (\ref{eq63}) can be replaced by 0. Thus we have
\begin{equation}
\tfrac{1}{2}\alpha^{2}\arcsin\left(  {\frac{\zeta}{\alpha}}\right)  +\tfrac
{1}{2}\zeta\left(  {\alpha^{2}-\zeta^{2}}\right)  ^{1/2}=\sigma E\left(
{x;\sigma^{-1}}\right)  , \label{eq65}%
\end{equation}
for $0\leq x\leq\sigma$ ($0\leq\zeta\leq\alpha$).

For $\sigma\leq x\leq1-\delta_{0}$ we have
\begin{equation}
\int_{\alpha}^{\zeta}{\left(  {\tau^{2}-\alpha^{2}}\right)  ^{1/2}d\tau}%
=\int_{\sigma}^{x}{\left\{  {f\left(  {\sigma,t}\right)  }\right\}  ^{1/2}%
dt}=\int_{\sigma}^{x}{\left(  {\frac{t^{2}-\sigma^{2}}{1-t^{2}}}\right)
^{1/2}dt}. \label{eq66}%
\end{equation}
Thus in this case
\begin{equation}
-\tfrac{1}{2}\alpha^{2}\operatorname{arccosh}\left(  {\frac{\zeta}{\alpha}%
}\right)  +\tfrac{1}{2}\zeta\left(  {\zeta^{2}-\alpha^{2}}\right)
^{1/2}=\left\vert \operatorname{Im}{\left\{  {\sigma E\left(  {x;\sigma^{-1}%
}\right)  }\right\}  }\right\vert . \label{eq67}%
\end{equation}

With
\begin{equation}
W=\left\{  {\frac{\sigma^{2}-x^{2}}{\left(  {\alpha^{2}-\zeta^{2}}\right)
\left(  {1-x^{2}}\right)  }}\right\}  ^{1/4}w, \label{eq68}%
\end{equation}
we transform (\ref{eq30}) to the form
\begin{equation}
\frac{d^{2}W}{d\zeta^{2}}=\left\{  {\gamma^{2}\left(  {\zeta^{2}-\alpha^{2}%
}\right)  +\psi\left(  {\gamma,\alpha,\zeta}\right)  }\right\}  W,
\label{eq69}%
\end{equation}
where
\begin{equation}%
\begin{array}
[c]{l}%
\psi\left(  {\gamma,\alpha,\zeta}\right)  =\dfrac{\left(  {1-m^{2}}\right)
\left(  {\alpha^{2}-\zeta^{2}}\right)  }{\left(  {1-x^{2}}\right)  \left(
{\sigma^{2}-x^{2}}\right)  }+\dfrac{2\alpha^{2}+3\zeta^{2}}{4\left(
{\alpha^{2}-\zeta^{2}}\right)  ^{2}}\\
-\dfrac{\left(  {1-\sigma^{2}}\right)  \left(  {\alpha^{2}-\zeta^{2}}\right)
\left\{  {6x^{4}-\left(  {\sigma^{2}+3}\right)  x^{2}-2\sigma^{2}}\right\}
}{4\left(  {1-x^{2}}\right)  \left(  {\sigma^{2}-x^{2}}\right)  ^{3}}.
\end{array}
\label{eq70}%
\end{equation}
To sharpen the subsequent error bounds it is possible to perturb the parameter
by defining a new parameter $\omega$ by $\alpha^{2}=\omega^{2}+\psi\left(
{\gamma,\alpha,0}\right)  \gamma^{-2}$, but we shall not pursue this here.

From theorem I of [26], with $u$ replaced by $\gamma$, we obtain two
independent solutions of (\ref{eq69}) given by
\begin{equation}
w_{1}\left(  {\gamma,\alpha,\zeta}\right)  =U\left(  {-{\tfrac{1}{2}}%
\gamma\alpha^{2},\zeta\sqrt{2\gamma}}\right)  +\varepsilon_{1}\left(
{\gamma,\alpha,\zeta}\right)  , \label{eq71}%
\end{equation}
and
\begin{equation}
w_{2}\left(  {\gamma,\alpha,\zeta}\right)  =\bar{{U}}\left(  {-{\tfrac{1}{2}%
}\gamma\alpha^{2},\zeta\sqrt{2\gamma}}\right)  +\varepsilon_{2}\left(
{\gamma,\alpha,\zeta}\right)  . \label{eq72}%
\end{equation}
Here $U\left(  {a,x}\right)  $ and $\bar{{U}}\left(  {a,x}\right)  $ are the
parabolic cylinder functions defined in [26, \S 5] and [27, \S 12.2], and are
linearly independent for $a<0$. The approximants $U\left(  {-{\frac{1}{2}%
}\gamma\alpha^{2},\zeta\sqrt{2\gamma}}\right)  $ and $\bar{{U}}\left(
{-{\frac{1}{2}}\gamma\alpha^{2},\zeta\sqrt{2\gamma}}\right)  $ satisfy the
differential equation (\ref{eq69}) with $\psi\left(  {\gamma,\alpha,\zeta
}\right)  \equiv0$.

The error terms are bounded by [26, \S 6], and in particular these show that
\begin{equation}
\varepsilon_{1}\left(  {\gamma,\alpha,\zeta}\right)  ={O}\left(
{\gamma^{-2/3}\ln\left(  \gamma\right)  }\right)  {\operatorname{env}}U\left(
{-{\tfrac{1}{2}}\gamma\alpha^{2},\zeta\sqrt{2\gamma}}\right)  , \label{eq73}%
\end{equation}
and
\begin{equation}
\varepsilon_{2}\left(  {\gamma,\alpha,\zeta}\right)  ={O}\left(
{\gamma^{-2/3}\ln\left(  \gamma\right)  }\right)  {\operatorname{env}}\bar
{{U}}\left(  {-{\tfrac{1}{2}}\gamma\alpha^{2},\zeta\sqrt{2\gamma}}\right)  ,
\label{eq74}%
\end{equation}
uniformly for $0\leq x\leq1-\delta_{0}$. Here the envelope function
${\operatorname{env}}$ is defined for the parabolic cylinder functions by [27,
eq. 14.15.23].

The parabolic cylinder function $U$ has the unique recessive property
\begin{equation}
U\left(  {-{\tfrac{1}{2}}a,x}\right)  \sim x^{\left(  {a-1}\right)
/2}e^{-x^{2}/4}\quad\left(  {x\rightarrow\infty}\right)  , \label{eq75}%
\end{equation}
whereas $\bar{{U}}$ is dominant, with the behavior
\begin{equation}
\bar{{U}}\left(  {-{\tfrac{1}{2}}a,x}\right)  \sim\left(  {2/\pi}\right)
^{1/2}\Gamma\left(  {{\tfrac{1}{2}}a+{\tfrac{1}{2}}}\right)  x^{-\left(
{a+1}\right)  /2}e^{x^{2}/4}\quad\left(  {x\rightarrow\infty}\right)  ;
\label{eq76}%
\end{equation}
see [26, \S 5]. In addition, from [26, Eqs. (5.12) and (5.13)]
and the definitions \eqref{eq71} and \eqref{eq72}, we note that $w_{1}\left(
{\gamma,\alpha,\zeta}\right)  $ and $w_{2}\left(  {\gamma,\alpha,\zeta
}\right)  $ are oscillatory in the $\zeta$ interval $\left[  {0,\alpha
}\right]  $, with comparable amplitudes and complementary phases of the argument.

Now, for negative $x$ and $\zeta,$ we will also need the solution given in
[26]
\begin{equation}
w_{4}\left(  {\gamma,\alpha,\zeta}\right)  =\bar{{U}}\left(  {-{\tfrac{1}{2}%
}\gamma\alpha^{2},-\zeta\sqrt{2\gamma}}\right)  +\varepsilon_{4}\left(
{\gamma,\alpha,\zeta}\right)  . \label{eq77}%
\end{equation}
We remark that
\begin{equation}
\varepsilon_{j}\left(  {\gamma,\alpha,0}\right)  =\partial\varepsilon
_{j}\left(  {\gamma,\alpha,0}\right)  /\partial\zeta=0\quad\left(
{j=2,4}\right)  , \label{eq78}%
\end{equation}
and hence
\begin{equation}
w_{2}\left(  {\gamma,\alpha,0}\right)  =w_{4}\left(  {\gamma,\alpha,0}\right)
=\bar{{U}}\left(  {-{\tfrac{1}{2}}\gamma\alpha^{2},0}\right)  , \label{eq79}%
\end{equation}
as well as
\begin{equation}
\partial w_{2}\left(  {\gamma,\alpha,0}\right)  /\partial\zeta=-\partial
w_{4}\left(  {\gamma,\alpha,0}\right)  /\partial\zeta=\sqrt{2\gamma}\bar
{{{U}^{\prime}}}\left(  {-{\tfrac{1}{2}}\gamma\alpha^{2},0}\right)  .
\label{eq80}%
\end{equation}
The error bounds for $\varepsilon_{4}\left(  {\gamma,\alpha,\zeta}\right)  $
only apply for non-positive $\zeta.$ In order to extend the solution to
positive values of $\zeta$ we use [26, Eqs. (6.23) and (6.24)] to obtain the
connection formula
\begin{equation}%
\begin{array}
[c]{l}%
w_{4}\left(  {\gamma,\alpha,\zeta}\right)  =-\left\{  {\sin\left(  {{\frac
{1}{2}}\pi\gamma\alpha^{2}}\right)  +{O}\left(  {\gamma^{-2/3}%
}\right)  }\right\}  w_{2}\left(  {\gamma,\alpha,\zeta}\right) \\
+\left\{  {\cos\left(  {{\frac{1}{2}}\pi\gamma\alpha^{2}}\right)
+{O}\left(  {\gamma^{-2/3}}\right)  }\right\}  w_{1}\left(
{\gamma,\alpha,\zeta}\right)  .
\end{array}
\label{eq81}%
\end{equation}

Now from (\ref{eq42}), (\ref{eq56}) and (\ref{eq64})
\begin{equation}
{\tfrac{1}{2}}\pi\gamma\alpha^{2}=\left(  {n-m+{\tfrac{1}{2}}}\right)
\pi+{O}\left(  {\gamma^{-1}}\right)  . \label{eq82}%
\end{equation}
Bearing in mind that $w_{1}\left(  {\gamma,\alpha,\zeta}\right)  $ is
exponentially small compared to $w_{2}\left(  {\gamma,\alpha,\zeta}\right)  $
in $0\leq x\leq1-\delta_{0}$ (except near its zeros) we deduce from
\eqref{eq81} and \eqref{eq82} that
\begin{equation}%
\begin{array}
[c]{l}%
w_{2}\left(  {\gamma,\alpha,\zeta}\right)  -\left(  {-1}\right)  ^{m+n}%
w_{4}\left(  {\gamma,\alpha,\zeta}\right)  =2w_{2}\left(  {\gamma,\alpha
,\zeta}\right) \\
+{O}\left(  {\gamma^{-2/3}}\right)  \left\{  {w_{1}\left(
{\gamma,\alpha,\zeta}\right)  +w_{2}\left(  {\gamma,\alpha,\zeta}\right)
}\right\}  .
\end{array}
\label{eq83}%
\end{equation}
We next express
\begin{equation}%
\begin{array}
[c]{l}%
\operatorname{Ps}_{n}^{m}\left(  {x,\gamma^{2}}\right)  =\left\{
{\dfrac{\alpha^{2}-\zeta^{2}}{\left(  {\sigma^{2}-x^{2}}\right)  \left(
{1-x^{2}}\right)  }}\right\}  ^{1/4}\\
\times\left[  {d_{n}^{m}\left(  \gamma\right)  w_{1}\left(  {\gamma
,\alpha,\zeta}\right)  +e_{n}^{m}\left(  \gamma\right)  \left\{  {w_{2}\left(
{\gamma,\alpha,\zeta}\right)  -\left(  {-1}\right)  ^{m+n}w_{4}\left(
{\gamma,\alpha,\zeta}\right)  }\right\}  }\right]  ,
\end{array}
\label{eq84}%
\end{equation}
and we shall determine the constant $d_{n}^{m}\left(  \gamma\right)  $ (as
well as bounding $e_{n}^{m}\left(  \gamma\right)  )$ by comparing both sides
of this relationship at appropriate values of $x$.

To this end, firstly we assume that $\operatorname{Ps}_{n}^{m}\left(
{x,\gamma^{2}}\right)  $ is even, so that $m+n$ is also even. Then, setting
$x=\zeta=0$ in (\ref{eq84}), and invoking (\ref{eq79}), immediately yields
\begin{equation}
d_{n}^{m}\left(  \gamma\right)  =\left(  {\frac{\sigma}{\alpha}}\right)
^{1/2}\frac{\operatorname{Ps}_{n}^{m}\left(  {0,\gamma^{2}}\right)  }%
{w_{1}\left(  {\gamma,\alpha,0}\right)  }. \label{eq85}%
\end{equation}
An asymptotic approximation for this constant, which does not involve
$\operatorname{Ps}_{n}^{m}\left(  {0,\gamma^{2}}\right)  $, is given by
(\ref{eq102}) below.

Next, if we differentiate both sides of (\ref{eq84}) and again set $x=\zeta=0$
we find from the property $\operatorname{Ps}_{n}^{m}{}^{\prime}\left(
{0,\gamma^{2}}\right)  =0$ that $e_{n}^{m}\left(  \gamma\right)
={O}\left(  {\gamma^{-1}}\right)  $, which is not sharp enough.
Instead we match the parabolic cylinder and Bessel function approximations,
and their derivatives, at the fixed point $x=1-{\frac{1}{2}}\delta_{0}$ (at
which both the parabolic cylinder function and modified Bessel function
approximations are valid). Using (\ref{eq61}) (\ref{eq83}) and (\ref{eq84}),
we therefore arrive at
\begin{equation}
e_{n}^{m}\left(  \gamma\right)  \sim-d_{n}^{m}\left(  \gamma\right)
\frac{\mathcal{W}\left\{  {\left(  {\zeta^{2}-\alpha^{2}}\right)
^{1/4}U\left(  {-{\frac{1}{2}}\gamma\alpha^{2},\zeta\sqrt{2\gamma}}\right)
,\left\vert \eta\right\vert ^{1/4}I_{m}\left(  {\gamma\left\vert
\eta\right\vert ^{1/2}}\right)  }\right\}  }{\mathcal{W}\left\{  {\left(
{\zeta^{2}-\alpha^{2}}\right)  ^{1/4}\overline{U}\left(  {-{\frac{1}{2}}%
\gamma\alpha^{2},\zeta\sqrt{2\gamma}}\right)  ,\left\vert \eta\right\vert
^{1/4}I_{m}\left(  {\gamma\left\vert \eta\right\vert ^{1/2}}\right)
}\right\}  }, \label{eq86}%
\end{equation}
and
\begin{equation}
c_{n}^{m}\left(  \gamma\right)  \sim d_{n}^{m}\left(  \gamma\right)
\frac{{\left(  {\zeta^{2}-\alpha^{2}}\right)  ^{1/2}}\mathcal{W}\left\{
{\overline{U}\left(  {-{\frac{1}{2}}\gamma\alpha^{2},\zeta\sqrt{2\gamma}%
}\right)  ,U\left(  {-{\frac{1}{2}}\gamma\alpha^{2},\zeta\sqrt{2\gamma}%
}\right)  }\right\}  }{\mathcal{W}\left\{  {\left(  {\zeta^{2}-\alpha^{2}%
}\right)  ^{1/4}\overline{U}\left(  {-{\frac{1}{2}}\gamma\alpha^{2},\zeta
\sqrt{2\gamma}}\right)  ,\left\vert \eta\right\vert ^{1/4}I_{m}\left(
{\gamma\left\vert \eta\right\vert ^{1/2}}\right)  }\right\}  }. \label{eq87}%
\end{equation}
In both of these the Wronskians $\mathcal{W}$ are with respect to $x$, and
evaluated at $x=1-{\frac{1}{2}}\delta_{0}$ (with $\eta$ and $\zeta$
corresponding to this value).

Next, from (\ref{eq66}) and [27, Eqs. 12.10.3 - 12.10.6], we have the
asymptotic approximations for large $\gamma$, fixed $\zeta\in\left(
{\alpha,\infty}\right)  $ and fixed $\alpha>0$
\begin{equation}
U\left(  {-{\tfrac{1}{2}}\gamma\alpha^{2},\zeta\sqrt{2\gamma}}\right)
\sim\left(  {\frac{\gamma\alpha^{2}}{2e}}\right)  ^{\gamma\alpha^{2}/4}%
\frac{\exp\left\{  {-\gamma\int_{\sigma}^{x}{\left\{  {f\left(  {\sigma
,t}\right)  }\right\}  ^{1/2}dt}}\right\}  }{\left\{  {2\gamma\left(
{\zeta^{2}-\alpha^{2}}\right)  }\right\}  ^{1/4}}, \label{eq88}%
\end{equation}%
\begin{equation}%
\begin{array}
[c]{l}%
{U}^{\prime}\left(  {-{\tfrac{1}{2}}\gamma\alpha^{2},\zeta\sqrt{2\gamma}%
}\right)  \sim-\frac{1}{2}\left(  {\dfrac{\gamma\alpha^{2}}{2e}}\right)
^{\gamma\alpha^{2}/4}\\
\times\left\{  {2\gamma\left(  {\zeta^{2}-\alpha^{2}}\right)  }\right\}
^{1/4}\exp\left\{  {-\gamma\int_{\sigma}^{x}{\left\{  {f\left(  {\sigma
,t}\right)  }\right\}  ^{1/2}dt}}\right\}  ,
\end{array}
\label{eq89}%
\end{equation}%
\begin{equation}
\overline{U}\left(  {-{\tfrac{1}{2}}\gamma\alpha^{2},\zeta\sqrt{2\gamma}%
}\right)  \sim2\left(  {\frac{\gamma\alpha^{2}}{2e}}\right)  ^{\gamma
\alpha^{2}/4}\frac{\exp\left\{  {\gamma\int_{\sigma}^{x}{\left\{  {f\left(
{\sigma,t}\right)  }\right\}  ^{1/2}dt}}\right\}  }{\left\{  {2\gamma\left(
{\zeta^{2}-\alpha^{2}}\right)  }\right\}  ^{1/4}}, \label{eq90}%
\end{equation}
and
\begin{equation}%
\begin{array}
[c]{l}%
\overline{U}^{\prime}\left(  {-{\tfrac{1}{2}}\gamma\alpha^{2},\zeta
\sqrt{2\gamma}}\right)  \sim\left(  {\dfrac{\gamma\alpha^{2}}{2e}}\right)
^{\gamma\alpha^{2}/4}\left\{  {2\gamma\left(  {\zeta^{2}-\alpha^{2}}\right)
}\right\}  ^{1/4}\\
\times\exp\left\{  {\gamma\int_{\sigma}^{x}{\left\{  {f\left(  {\sigma
,t}\right)  }\right\}  ^{1/2}dt}}\right\}  .
\end{array}
\label{eq91}%
\end{equation}

These, along with
\begin{equation}
\left\vert \eta\right\vert ^{1/4}I_{m}\left(  {\gamma\left\vert \eta
\right\vert ^{1/2}}\right)  \sim\left(  {2\pi\gamma}\right)  ^{-1/2}%
\exp\left\{  {\gamma\int_{x}^{1}{\left\{  {f\left(  {\sigma,t}\right)
}\right\}  ^{1/2}dt}}\right\}  , \label{eq92}%
\end{equation}%
\begin{equation}
\frac{d\left\{  {\left\vert \eta\right\vert ^{1/4}I_{m}\left(  {\gamma
\left\vert \eta\right\vert ^{1/2}}\right)  }\right\}  }{dx}\sim-\left(
{\frac{\gamma}{2\pi}}\right)  ^{1/2}\left(  {\frac{x^{2}-\sigma^{2}}{1-x^{2}}%
}\right)  ^{1/2}\exp\left\{  {\gamma\int_{x}^{1}{\left\{  {f\left(  {\sigma
,t}\right)  }\right\}  ^{1/2}dt}}\right\}  , \label{eq93}%
\end{equation}
and
\begin{equation}
\frac{d\zeta}{dx}=\left\{  {\frac{x^{2}-\sigma^{2}}{\left(  {1-x^{2}}\right)
\left(  {\zeta^{2}-\alpha^{2}}\right)  }}\right\}  ^{1/2}, \label{eq94}%
\end{equation}
can be used to simplify (\ref{eq86}) and (\ref{eq87}). In particular, we find
that
\begin{equation}%
\begin{array}
[c]{l}%
e_{n}^{m}\left(  \gamma\right)  \left\{  {w_{2}\left(  {\gamma,\alpha,\zeta
}\right)  -\left(  {-1}\right)  ^{m+n}w_{4}\left(  {\gamma,\alpha,\zeta
}\right)  }\right\} \\
=o\left(  1\right)  A{\operatorname{env}}U\left(  {-{\tfrac{1}{2}}\gamma
\alpha^{2},\zeta\sqrt{2\gamma}}\right)  ,
\end{array}
\label{eq95}%
\end{equation}
where the $o\left(  1\right)  $ term is exponentially small as $\gamma
\rightarrow\infty_{\,}$for $x\in\left[  {0,1-\delta_{0}}\right]  $. In
addition, we obtain the useful result
\begin{equation}
c_{n}^{m}\left(  \gamma\right)  \sim d_{n}^{m}\left(  \gamma\right)  \left(
{\frac{\gamma\alpha^{2}}{2e}}\right)  ^{\gamma\alpha^{2}/4}\left(  {\frac
{\pi^{2}}{2\gamma}}\right)  ^{1/4}\exp\left\{  {-\gamma\int_{\sigma}%
^{1}{\left\{  {f\left(  {\sigma,t}\right)  }\right\}  ^{1/2}dt}}\right\}  .
\label{eq96}%
\end{equation}
From (\ref{eq73}) and (\ref{eq84}) - (\ref{eq95}), for $m+n$ even, $m$ bounded
and $n$ satisfying (\ref{eq2}), we arrive at our desired result
\begin{equation}%
\begin{array}
[c]{l}%
\operatorname{Ps}_{n}^{m}\left(  {x,\gamma^{2}}\right)  =\dfrac
{\operatorname{Ps}\left(  {0,\gamma^{2}}\right)  }{U\left(  {-{\frac{1}{2}%
}\gamma\alpha^{2},0}\right)  }\left\{  {\dfrac{\sigma^{2}\left(  {\alpha
^{2}-\zeta^{2}}\right)  }{\alpha^{2}\left(  {\sigma^{2}-x^{2}}\right)  \left(
{1-x^{2}}\right)  }}\right\}  ^{1/4}\\
\times\left\{  {U\left(  {-{\frac{1}{2}}\gamma\alpha^{2},\zeta\sqrt{2\gamma}%
}\right)  +{O}\left(  {\gamma^{-2/3}\ln\left(  \gamma\right)
}\right)  \operatorname{env}U\left(  {-{\frac{1}{2}}\gamma\alpha^{2}%
,\zeta\sqrt{2\gamma}}\right)  }\right\}  ,
\end{array}
\label{eq97}%
\end{equation}
as $\gamma\rightarrow\infty$, uniformly for $0\leq x\leq1-\delta_{0}$.

From [26, \S 5] we note that
\begin{equation}
U\left(  {-{\tfrac{1}{2}}\gamma\alpha^{2},0}\right)  =\pi^{-1/2}2^{\left(
{\gamma\alpha^{2}-1}\right)  /4}\Gamma\left(  {{\frac{1}{4}}\gamma\alpha
^{2}+{\frac{1}{4}}}\right)  \sin\left(  {{\frac{1}{4}}\gamma\alpha^{2}%
\pi+{\frac{1}{4}}\pi}\right)  , \label{eq98}%
\end{equation}
as well as
\begin{equation}
{U}^{\prime}\left(  {-{\tfrac{1}{2}}\gamma\alpha^{2},0}\right)  =-\pi
^{-1/2}2^{\left(  {\gamma\alpha^{2}+1}\right)  /4}\Gamma\left(  {{\frac{1}{4}%
}\gamma\alpha^{2}+{\frac{3}{4}}}\right)  \sin\left(  {{\frac{1}{4}}%
\gamma\alpha^{2}\pi+{\frac{3}{4}}\pi}\right)  . \label{eq99}%
\end{equation}
Thus, on referring to (\ref{eq82}), we observe that the RHS of (\ref{eq98}) is
bounded away from zero for large $\gamma$ when $m+n$ is even, and likewise for
the RHS of \eqref{eq99} when $m+n$ is odd (see \eqref{eq101} below).

For the case $\operatorname{Ps}_{n}^{m}\left(  {x,\gamma^{2}}\right)  $ odd,
equivalently $m+n$ odd, we differentiate both sides of \eqref{eq84} with
respect to $\zeta$, and then set $x=\zeta=0$. As a result, using \eqref{eq62}
and \eqref{eq80}, along with the fact that $\operatorname{Ps}_{n}^{m}\left(
{0,\gamma^{2}}\right)  =0$, we obtain
\begin{equation}
d_{n}^{m}\left(  \gamma\right)  =\left(  {\frac{\alpha}{\sigma}}\right)
^{1/2}\frac{\operatorname{Ps}_{n}^{m}{}^{\prime}\left(  {0,\gamma^{2}}\right)
}{\partial w_{1}\left(  {\gamma,\alpha,0}\right)  /\partial\zeta}.
\label{eq100}%
\end{equation}
Thus, again from \eqref{eq95}, we conclude for $m+n$ odd, $m$ bounded and $n$
satisfying \eqref{eq2}, that
\begin{equation}%
\begin{array}
[c]{l}%
\operatorname{Ps}_{n}^{m}\left(  {x,\gamma^{2}}\right)  =\dfrac
{\operatorname{Ps}_{n}^{m}{}^{\prime}\left(  {0,\gamma^{2}}\right)  }%
{{U}^{\prime}\left(  {-{\frac{1}{2}}\gamma\alpha^{2},0}\right)  }\left\{
{\dfrac{\alpha^{2}\left(  {\alpha^{2}-\zeta^{2}}\right)  }{4\gamma^{2}%
\sigma^{2}\left(  {\sigma^{2}-x^{2}}\right)  \left(  {1-x^{2}}\right)  }%
}\right\}  ^{1/4}\\
\times\left\{  {U\left(  {-{\frac{1}{2}}\gamma\alpha^{2},\zeta\sqrt{2\gamma}%
}\right)  +{O}\left(  {\gamma^{-2/3}\ln\left(  \gamma\right)
}\right)  \operatorname{env}U\left(  {-{\frac{1}{2}}\gamma\alpha^{2}%
,\zeta\sqrt{2\gamma}}\right)  }\right\}  ,
\end{array}
\label{eq101}%
\end{equation}
as $\gamma\rightarrow\infty$, uniformly for $0\leq x\leq1-\delta$. In this
${U}^{\prime}\left(  {-{\frac{1}{2}}\gamma\alpha^{2},0}\right)  $ is given by
(\ref{eq99}).

We now show that the proportionality constants in (\ref{eq97}) and
(\ref{eq101}) can be replaced by one that does not involve $\operatorname{Ps}%
_{n}^{m}\left(  {0,\gamma^{2}}\right)  $ or $\operatorname{Ps}_{n}^{m}%
{}^{\prime}\left(  {0,\gamma^{2}}\right)  $. Specifically, from (\ref{eq19}),
(\ref{eq61}), (\ref{eq84}), (\ref{eq95}) and (\ref{eq96}) we have (for both
the even and odd cases) that
\begin{equation}
d_{n}^{m}\left(  \gamma\right)  \sim\left\{  {\frac{\left(  {n+m}\right)
!}{\left(  {2n+1}\right)  \left(  {n-m}\right)  !p_{n}^{m}\left(
\gamma\right)  }}\right\}  ^{1/2}, \label{eq102}%
\end{equation}
as $\gamma\rightarrow\infty$, again with $m$ bounded and $n$ satisfying
(\ref{eq2}). Here
\begin{equation}%
\begin{array}
[c]{l}%
p_{n}^{m}\left(  \gamma\right)  =\left[  {\int_{0}^{1-\delta_{0}}{\left\{
{\dfrac{\alpha^{2}-\zeta^{2}}{\left(  {\sigma^{2}-x^{2}}\right)  \left(
{1-x^{2}}\right)  }}\right\}  ^{1/2}U^{2}\left(  {-{\frac{1}{2}}\gamma
\alpha^{2},\zeta\sqrt{2\gamma}}\right)  dx}}\right. \\
\left.  {+q_{n}^{m}\left(  \gamma\right)  \int_{1-\delta_{0}}^{1}{\left\{
{\dfrac{\left\vert \eta\right\vert }{\left(  {1-x^{2}}\right)  \left(
{x^{2}-\sigma^{2}}\right)  }}\right\}  ^{1/2}I_{m}^{2}\left(  {\gamma
\left\vert \eta\right\vert ^{1/2}}\right)  dx}}\right]  ,
\end{array}
\label{eq103}%
\end{equation}
in which
\begin{equation}
q_{n}^{m}\left(  \gamma\right)  =\left(  {\frac{\gamma\alpha^{2}}{2e}}\right)
^{\gamma\alpha^{2}/2}\left(  {\frac{\pi^{2}}{2\gamma}}\right)  ^{1/2}%
\exp\left\{  {-2\gamma\int_{\sigma}^{1}{\left\{  {f\left(  {\sigma,t}\right)
}\right\}  ^{1/2}dt}}\right\}  . \label{eq104}%
\end{equation}
Note also, from (\ref{eq96}), that under the same conditions
\begin{equation}
c_{n}^{m}\left(  \gamma\right)  \sim\left\{  {\frac{\left(  {n+m}\right)
!q_{n}^{m}\left(  \gamma\right)  }{\left(  {2n+1}\right)  \left(
{n-m}\right)  !p_{n}^{m}\left(  \gamma\right)  }}\right\}  ^{1/2}.
\label{eq105}%
\end{equation}

\section{Fixed $m$ and $n$: the angular case}

For fixed $m$\textit{ and }$n$ we can simplify the results of the previous
section, by applying the theory of [10]. To this end we observe that
(\ref{eq26}) can be expressed in the form
\begin{equation}
\frac{d^{2}w}{dx^{2}}=\left[  {\frac{\gamma^{2}x^{2}}{1-x^{2}}-\frac{a\gamma
}{1-x^{2}}+\frac{m^{2}-1}{\left(  {1-x^{2}}\right)  ^{2}}}\right]  w,
\label{eq106}%
\end{equation}
where
\begin{equation}
a=\lambda\gamma^{-1}+\gamma=2\left(  {n-m+{\tfrac{1}{2}}}\right)
+{O}\left(  {\gamma^{-1}}\right)  , \label{eq107}%
\end{equation}
the ${O}\left(  {\gamma^{-1}}\right)  $ term being valid for fixed $m$
and $n$ and $\gamma\rightarrow\infty$. In particular, $a$ is bounded.

Equation (\ref{eq106}) is characterised as having a pair of almost coalescent
turning points near $x=0$. The appropriate Liouville transformation in this
case is given by
\begin{equation}
\frac{1}{2}\rho^{2}=\int_{0}^{x}{\frac{t}{\left(  {1-t^{2}}\right)  ^{1/2}}%
dt}=1-\left(  {1-x^{2}}\right)  ^{1/2}. \label{eq108}%
\end{equation}
Note $x=0$ corresponds to $\rho=0$, and $x=1$ corresponds to $\rho=\sqrt{2}$.
Then with
\begin{equation}
W=\frac{x^{1/2}}{\rho^{1/2}\left(  {1-x^{2}}\right)  ^{1/4}}w, \label{eq109}%
\end{equation}
we obtain
\begin{equation}
\frac{d^{2}W}{d\rho^{2}}=\left[  {\gamma^{2}\rho^{2}-\gamma a+\gamma\zeta
\phi\left(  \rho\right)  +\chi\left(  \rho\right)  }\right]  W, \label{eq110}%
\end{equation}
where
\begin{equation}
\phi\left(  \rho\right)  =-\frac{a\rho}{4-\zeta^{2}}, \label{eq111}%
\end{equation}
and
\begin{equation}
\chi\left(  \rho\right)  =\frac{\rho^{2}\left(  {4m^{2}-1}\right)  }{\left(
{2-\rho^{2}}\right)  ^{2}}+\frac{7\rho^{2}-40}{4\left(  {4-\rho^{2}}\right)
^{2}}+\frac{4m^{2}}{\left(  {4-\rho^{2}}\right)  }. \label{eq112}%
\end{equation}
We remark that $\chi\left(  \rho\right)  ={O}\left(  1\right)  _{\,}%
$as $\gamma\rightarrow\infty$, and this function is analytic at $\rho=0$
($x=0$), but is not analytic at $\rho=\sqrt{2}\ $(${x=1}$).

Our approximants are again the parabolic cylinder functions $U\left(
{-{\frac{1}{2}}a,\rho\sqrt{2\gamma}}\right)  $ and $\bar{{U}}\left(
{-{\frac{1}{2}}a,\rho\sqrt{2\gamma}}\right)  $ (c.f. (\ref{eq71}) and
(\ref{eq72})). In this form they are solutions of
\begin{equation}
\frac{d^{2}W}{d\rho^{2}}=\left[  {\gamma^{2}\rho^{2}-\gamma a}\right]  W.
\label{eq113}%
\end{equation}
On comparing this equation with (\ref{eq110}) we note the extra
\textquotedblleft large\textquotedblright\ term $\gamma\zeta\phi\left(
\rho\right)  $. On account of this discrepancy we perturb the independent
variable, thus taking as approximants
\begin{equation}
U_{1}=\left\{  {1+\gamma^{-1}{\Phi}^{\prime}\left(  \rho\right)  }\right\}
^{-1/2}U\left(  {-{\tfrac{1}{2}}a,\hat{{\rho}}\sqrt{2\gamma}}\right)  ,
\label{eq114}%
\end{equation}
and
\begin{equation}
U_{2}=\left\{  {1+\gamma^{-1}{\Phi}^{\prime}\left(  \rho\right)  }\right\}
^{-1/2}\bar{{U}}\left(  {-{\tfrac{1}{2}}a,\hat{{\rho}}\sqrt{2\gamma}}\right)
, \label{eq115}%
\end{equation}
where
\begin{equation}
\hat{{\rho}}=\rho+\gamma^{-1}\Phi\left(  \rho\right)  , \label{eq116}%
\end{equation}
in which
\begin{equation}
\Phi\left(  \rho\right)  =\frac{1}{2\rho}\int_{0}^{\rho}{\phi\left(  v\right)
dv}=\frac{a\ln\left(  {1-{\frac{1}{4}}\rho^{2}}\right)  }{4\rho}.
\label{eq117}%
\end{equation}
In [10] it is shown that $U_{j}$ satisfy the differential equation
\begin{equation}
\frac{d^{2}U}{d\rho^{2}}=\left\{  {\gamma^{2}\rho^{2}-\gamma a+\gamma\rho
\phi\left(  \rho\right)  +g\left(  {\gamma,\rho}\right)  }\right\}  U,
\label{eq118}%
\end{equation}
where $g\left(  {\gamma,\rho}\right)  ={O}\left(  1\right)  _{\,}$as
$\gamma\rightarrow\infty$, uniformly for $\rho\in\left[  {0,\sqrt{2}-\delta
}\right]  $. Thus (\ref{eq118}) is the appropriate comparison equation to
(\ref{eq110}).

Following [10] we then define
\begin{equation}
\hat{{w}}_{j}\left(  {\gamma,\rho}\right)  =U_{j}\left(  {\gamma,\rho}\right)
+\hat{{\varepsilon}}_{j}\left(  {\gamma,\rho}\right)  \quad\left(
{j=1,2}\right)  , \label{eq119}%
\end{equation}
as exact solutions of (\ref{eq110}). Explicit error bounds are furnished in
[10], and from these it follows that
\begin{equation}
\hat{{\varepsilon}}_{1}\left(  {\gamma,\rho}\right)  ={O}\left(
{\gamma^{-1}\ln\left(  \gamma\right)  }\right)  {\operatorname{env}}U\left(
{-{\tfrac{1}{2}}a,\hat{{\rho}}\sqrt{2\gamma}}\right)  , \label{eq120}%
\end{equation}
uniformly for $0\leq x\leq1-\delta_{0}$, and similarly for $\hat{{\varepsilon
}}_{2}\left(  {\gamma,\rho}\right)  $.

Let us assume that $\operatorname{Ps}_{n}^{m}\left(  {x,\gamma^{2}}\right)  $
(and hence $m+n)$ is even. Similarly to (\ref{eq84}) we write
\begin{equation}%
\begin{array}
[c]{l}%
\operatorname{Ps}_{n}^{m}\left(  {x,\gamma^{2}}\right)  =\rho^{1/2}%
x^{-1/2}\left(  {1-x^{2}}\right)  ^{-1/4}\\
\times\left[  {\hat{{d}}_{n}^{m}\left(  \gamma\right)  \hat{{w}}_{1}\left(
{\gamma,\rho}\right)  +\hat{{e}}_{n}^{m}\left(  \gamma\right)  \left\{
{\hat{{w}}_{2}\left(  {\gamma,\rho}\right)  -\hat{{w}}_{4}\left(  {\gamma
,\rho}\right)  }\right\}  }\right]  ,
\end{array}
\label{eq121}%
\end{equation}
where $\hat{{w}}_{4}\left(  {\gamma,\rho}\right)  $ is the solution (involving
$\bar{{U}})$ given by eq. (110) of [10]. By matching at $x=\rho=0$ we find
\begin{equation}
\hat{{d}}_{n}^{m}\left(  \gamma\right)  =\frac{Ps_{n}^{m}\left(  {0,\gamma
^{2}}\right)  }{\hat{{w}}_{1}\left(  {\gamma,0}\right)  }. \label{eq122}%
\end{equation}
Analogously to the proof of (\ref{eq95}) it can be shown that
\begin{equation}
\hat{{e}}_{n}^{m}\left(  \gamma\right)  \left\{  {\hat{{w}}_{2}\left(
{\gamma,\rho}\right)  -\hat{{w}}_{4}\left(  {\gamma,\rho}\right)  }\right\}
=o\left(  1\right)  \hat{{A}}{\operatorname{env}}U\left(  {-{\tfrac{1}{2}%
}a,\hat{{\rho}}\sqrt{2\gamma}}\right)  , \label{eq123}%
\end{equation}
where $o\left(  1\right)  $ is exponentially small for $0\leq x\leq
1-\delta_{0}$ as $\gamma\rightarrow\infty$. Consequently, we arrive at our
desired result
\begin{equation}%
\begin{array}
[c]{l}%
\operatorname{Ps}_{n}^{m}\left(  {x,\gamma^{2}}\right)  =\dfrac
{\operatorname{Ps}_{n}^{m}\left(  {0,\gamma^{2}}\right)  }{U\left(
{-{\frac{1}{2}}a,0}\right)  }\left(  {\dfrac{\rho}{x}}\right)  ^{1/2}\left(
{1-x^{2}}\right)  ^{-1/4}\\
\times\left[  {U\left(  {-{\frac{1}{2}}a,\hat{{\rho}}\sqrt{2\gamma}}\right)
+{O}\left(  {\gamma^{-1}\ln\left(  \gamma\right)  }\right)
\operatorname{env}U\left(  {-{\frac{1}{2}}a,\hat{{\rho}}\sqrt{2\gamma}%
}\right)  }\right]  ,
\end{array}
\label{eq124}%
\end{equation}
as $\gamma\rightarrow\infty$, uniformly for $0\leq x\leq1-\delta_{0}$.

For the case $\operatorname{Ps}_{n}^{m}\left(  {x,\gamma^{2}}\right)  $ being
odd we likewise obtain, under the same conditions,
\begin{equation}%
\begin{array}
[c]{l}%
\operatorname{Ps}_{n}^{m}\left(  {x,\gamma^{2}}\right)  =\dfrac
{\operatorname{Ps}_{n}^{m}{}^{\prime}\left(  {0,\gamma^{2}}\right)  }%
{{U}^{\prime}\left(  {-{\frac{1}{2}}a,0}\right)  }\left(  {\dfrac{\rho
}{2\gamma x}}\right)  ^{1/2}\left(  {1-x^{2}}\right)  ^{-1/4}\\
\times\left[  {U\left(  {-{\frac{1}{2}}a,\hat{{\rho}}\sqrt{2\gamma}}\right)
+{O}\left(  {\gamma^{-1}\ln\left(  \gamma\right)  }\right)
\operatorname{env}U\left(  {-{\frac{1}{2}}a,\hat{{\rho}}\sqrt{2\gamma}%
}\right)  }\right]  .
\end{array}
\label{eq125}%
\end{equation}

\section{Summary}

For reference we collect the principal results of the paper. All results are
uniformly valid for $\gamma\to\infty$, $m$ and $n$ integers, $m$ bounded, and
$n$ satisfying $0\le m\le n\le2\pi^{-1}\gamma\left(  {1-\delta} \right)  $
where $\delta\in\left(  {0,1} \right)  $ is fixed.

We define $\sigma=\sqrt{1+\gamma^{-2}\lambda_{n}^{m}\left(  {\gamma^{2}%
}\right)  }$ and assume $0\leq\sigma\leq\sigma_{0}<1$ for an arbitrary fixed
positive $\sigma_{0}$. We further define variables $\xi=\xi\left(  x\right)  $
and $\zeta=\zeta\left(  x\right)  $ by
\begin{equation}
\xi=\int_{1}^{x}{\left(  {\frac{t^{2}-\sigma^{2}}{t^{2}-1}}\right)  ^{1/2}dt},
\label{eq126}%
\end{equation}
and
\begin{equation}
\int_{\alpha}^{\zeta}{\left\vert {\tau^{2}-\alpha^{2}}\right\vert ^{1/2}d\tau
}=\int_{\sigma}^{x}{\left(  {\frac{\left\vert {t^{2}-\sigma^{2}}\right\vert
}{1-t^{2}}}\right)  ^{1/2}dt}, \label{eq127}%
\end{equation}
where
\begin{equation}
\alpha=2\left\{  {\frac{1}{\pi}\int_{0}^{\sigma}{\left(  {\frac{\sigma
^{2}-t^{2}}{1-t^{2}}}\right)  ^{1/2}dt}}\right\}  ^{1/2}. \label{eq128}%
\end{equation}
Then, using the definition above for $\sigma$, a uniform asymptotic
relationship between $\lambda_{n}^{m}\left(  {\gamma^{2}}\right)  $ and the
parameters $m$, $n$ and $\gamma$ is given implicitly by the relation
\begin{equation}
\gamma\int_{0}^{\sigma}{\left(  {\frac{\sigma^{2}-t^{2}}{1-t^{2}}}\right)
^{1/2}dt}=\frac{1}{2}\left(  {n-m+\frac{1}{2}}\right)  \pi+{O}\left(
{\frac{1}{\gamma}}\right)  . \label{eq129}%
\end{equation}

The following approximation holds for the radial PSWF
\begin{equation}%
\begin{array}
[c]{l}%
Ps_{n}^{m}\left(  {x,\gamma^{2}}\right)  =\left\{  {\dfrac{\left(
{n+m}\right)  !q_{n}^{m}\left(  \gamma\right)  }{\left(  {2n+1}\right)
\left(  {n-m}\right)  !p_{n}^{m}\left(  \gamma\right)  }}\right\}
^{1/2}\left\{  {\left(  {x^{2}-1}\right)  \left(  {x^{2}-\sigma^{2}}\right)
}\right\}  ^{-1/4}\\
\times\xi^{1/2}\left[  {J_{m}\left(  {\gamma\xi}\right)  +{O}\left(
{\gamma^{-1}}\right)  \operatorname{env}J_{m}\left(  {\gamma\xi}\right)
}\right]  ,
\end{array}
\label{eq130}%
\end{equation}
this being uniformly valid for $1<x<\infty$. Here $J_{m}$ is the Bessel
function of the first kind, ${\operatorname{env}}J_{m}$ is defined by [27,
\S 2.8(iv)], and the constants $p_{n}^{m}\left(  \gamma\right)  $ and
$q_{n}^{m}\left(  \gamma\right)  $ are given by
\begin{equation}%
\begin{array}
[c]{l}%
p_{n}^{m}\left(  \gamma\right)  =\left[  {\int_{0}^{1-\delta_{0}}{\left\{
{\dfrac{\alpha^{2}-\zeta^{2}}{\left(  {\sigma^{2}-x^{2}}\right)  \left(
{1-x^{2}}\right)  }}\right\}  ^{1/2}U^{2}\left(  {-{\frac{1}{2}}\gamma
\alpha^{2},\zeta\sqrt{2\gamma}}\right)  dx}}\right. \\
\left.  {+q_{n}^{m}\left(  \gamma\right)  \int_{1-\delta_{0}}^{1}{\left\{
{\dfrac{1}{\left(  {1-x^{2}}\right)  \left(  {x^{2}-\sigma^{2}}\right)  }%
}\right\}  ^{1/2}\left\vert \xi\right\vert I_{m}^{2}\left(  {\gamma\left\vert
\xi\right\vert }\right)  dx}}\right]  ,
\end{array}
\label{eq131}%
\end{equation}
and
\begin{equation}
q_{n}^{m}\left(  \gamma\right)  =\left(  {\frac{\gamma\alpha^{2}}{2e}}\right)
^{\gamma\alpha^{2}/2}\left(  {\frac{\pi^{2}}{2\gamma}}\right)  ^{1/2}%
\exp\left\{  {-2\gamma\int_{\sigma}^{1}{\left(  {\frac{t^{2}-\sigma^{2}%
}{1-t^{2}}}\right)  ^{1/2}dt}}\right\}  . \label{eq132}%
\end{equation}
In (\ref{eq131}) $\delta_{0}\in\left(  {0,1-\sigma_{0}}\right)  $ is
arbitrarily chosen, $I_{m}$ is the modified Bessel function of the first kind,
and $U$ is the parabolic cylinder function (see [26, \S 5]).

In terms of the modified Bessel function, we have for the radial PSWF
\begin{equation}%
\begin{array}
[c]{l}%
\operatorname{Ps}_{n}^{m}\left(  {x,\gamma^{2}}\right)  =\left\{
{\dfrac{\left(  {n+m}\right)  !q_{n}^{m}\left(  \gamma\right)  }{\left(
{2n+1}\right)  \left(  {n-m}\right)  !p_{n}^{m}\left(  \gamma\right)  }%
}\right\}  ^{1/2}\left\{  {\left(  {1-x^{2}}\right)  \left(  {x^{2}-\sigma
^{2}}\right)  }\right\}  ^{-1/4}\\
\times\left\vert \xi\right\vert ^{1/2}I_{m}\left(  {\gamma\left\vert
\xi\right\vert }\right)  \left\{  {1+{O}\left(  {\gamma^{-1}}\right)
}\right\}  ,
\end{array}
\label{eq133}%
\end{equation}
uniformly for $1-\delta_{0}\leq x<1$.

Finally, in terms of the parabolic cylinder function, the asymptotic
approximation
\begin{equation}%
\begin{array}
[c]{l}%
\operatorname{Ps}_{n}^{m}\left(  {x,\gamma^{2}}\right)  =\left\{
{\dfrac{\left(  {n+m}\right)  !}{\left(  {2n+1}\right)  \left(  {n-m}\right)
!p_{n}^{m}\left(  \gamma\right)  }}\right\}  ^{1/2}\left\{  {\dfrac{\alpha
^{2}-\zeta^{2}}{\left(  {\sigma^{2}-x^{2}}\right)  \left(  {1-x^{2}}\right)
}}\right\}  ^{1/4}\\
\times\left\{  {U\left(  {-{\frac{1}{2}}\gamma\alpha^{2},\zeta\sqrt{2\gamma}%
}\right)  +{O}\left(  {\gamma^{-2/3}\ln\left(  \gamma\right)
}\right)  \operatorname{env}U\left(  {-{\frac{1}{2}}\gamma\alpha^{2}%
,\zeta\sqrt{2\gamma}}\right)  }\right\}  ,
\end{array}
\label{eq134}%
\end{equation}
holds uniformly for $0\leq x\leq1-\delta_{0}$, where ${\operatorname{env}}U$
is defined by [27, eq. 14.15.23].

\section*{Acknowledgement}
        I thank the referee for a number helpful suggestions.

%%%%%%%% Bibliography %%%%%%%%%%%%%%%%%%%%%%%%%%%%%%%%%%%%%%%%%%%%%%%%%

\end{document}